\documentclass[12pt]{amsart}
\usepackage{amssymb,latexsym,eufrak,amsmath,amscd, graphicx}
\usepackage[colorlinks=true,pagebackref,hyperindex]{hyperref}

\usepackage[all]{xy}
\usepackage{amsfonts}
\usepackage{amscd}

\usepackage{xypic}

\renewcommand{\phi}{\varphi}
\renewcommand{\epsilon}{\varepsilon}
\renewcommand{\theta}{\vartheta}

\def\ZZ{{\mathbf Z}}
\def\NN{{\mathbf N}}

\def\AAA{{\mathbf A}}
\def\RR{{\mathbf R}}
\def\QQ{{\mathbf Q}}
\def\PP{{\mathbf P}}

\newcommand{\leftidx}[3]{{\vphantom{#2}}#1#2#3}

\newcommand{\ltightstar}[1]{\leftidx{^*}{\negthinspace#1}{}}

\def\cJ{\mathcal{J}}

\def\cF{\mathcal{F}}

\def\cG{\mathcal{G}}
\def\cL{\mathcal{L}}
\def\cO{\mathcal{O}}
\def\cH{\mathcal{H}}
\def\cT{\mathcal{T}}

\def\cU{\mathcal{U}}

\def\fra{\mathfrak{a}}
\def\frb{\mathfrak{b}}

\def\frm{\mathfrak{m}}

 \DeclareMathOperator{\Spec}{Spec}
 
 \DeclareMathOperator{\lct}{lct}
 
 \DeclareMathOperator{\fpt}{fpt}

\newcommand{\llbracket}{[\negthinspace[}
\newcommand{\rrbracket}{]\negthinspace]}

\newtheorem{lemma}{Lemma}[section]
\newtheorem{theorem}[lemma]{Theorem}
\newtheorem{corollary}[lemma]{Corollary}
\newtheorem{proposition}[lemma]{Proposition}

\theoremstyle{definition}

\newtheorem{remark}[lemma]{Remark}

\theoremstyle{remark}
\newtheorem*{remark*}{Remark}
\newtheorem*{note*}{Note}

\begin{document}

\title{Log canonical thresholds, $F$-pure thresholds, and non-standard extensions}

\thanks{2000\,\emph{Mathematics Subject Classification}.
 Primary 13A35; Secondary 13L05, 14B05, 14F18.
\newline The second author was partially supported by
RTG grant 0502170, and the fourth author was partially supported by
 NSF grant DMS-0758454 and
  a Packard Fellowship.}

\author[B.~Bhatt]{Bhargav~Bhatt}
\address{Department of Mathematics, University of Michigan,
Ann Arbor, MI 48109, USA}
\email{{bhargav.bhatt@gmail.com}}

\author[D.~J.~Hern\'{a}ndez]{Daniel~J.~Hern\'{a}ndez}
\address{Department of Mathematics, University of Michigan,
Ann Arbor, MI 48109, USA}
\email{{dhernan@umich.edu}}

\author[L.~E.~Miller]{Lance~E.~Miller}
\address{Department of Mathematics, University of Utah, 155 South 1400 East,
Salt Lake City, UT 84112}
\email{{lmiller@math.utah.edu}}

\author[M.~Musta\c{t}\u{a}]{Mircea~Musta\c{t}\u{a}}
\address{Department of Mathematics, University of Michigan,
Ann Arbor, MI 48109, USA}
\email{{mmustata@umich.edu}}

\begin{abstract}
We present a new relation between an invariant of singularities in
characteristic zero (the log canonical threshold) and an invariant of singularities
defined via the Frobenius morphism in positive characteristic (the $F$-pure threshold). We show that the set of limit points 
of sequences of the form $(c_p)$, where $c_p$ is the $F$-pure threshold of an
ideal on an $n$-dimensional smooth variety
in characteristic $p$, coincides with the set of log canonical thresholds of ideals
on $n$-dimensional smooth varieties
in characteristic zero. We prove this by combining results of Hara and Yoshida
with non-standard constructions.
\end{abstract}

\maketitle

\markboth{B.~BHATT, D.~J.~HERN\'{A}NDEZ, L.~E.~MILLER, AND M.~MUSTA\c{T}\u{A}}
{LOG CANONICAL AND $F$-PURE THRESHOLDS, AND NON-STANDARD EXTENSIONS}

\section{Introduction}\label{intro}

The connection between invariants of singularities in characteristic zero and positive characteristic
is a topic that has recently attracted a lot of attention. Typically, the invariants of singularities that arise in birational geometry are defined via divisorial valuations. In characteristic zero, one can use
(log) resolutions of singularities to compute such invariants. On the other hand, in commutative algebra in positive characteristic one defines invariants using the action of the Frobenius morphism. 
It turns out that these invariants have subtle connections, some of them proven, and some still conjectural (see, for example, 
\cite{HW}, \cite{HY}, and \cite{MTW}). The typical such connection involves reduction from characteristic zero to positive characteristic. In this note we describe a different, though related connection. We use non-standard constructions to study
limits of invariants in positive characteristic, where the characteristic tends to infinity, in terms of 
invariants in characteristic zero.

The invariants we study in this paper are the log canonical threshold (in characteristic zero)
and the $F$-pure threshold (in positive characteristic). The log canonical threshold is an
invariant that plays an important role in birational geometry (see \cite{Kol} and 
\cite{EM}). Given an irreducible, smooth scheme $X$ defined over a field
$k$ of characteristic zero, and a proper ideal $\fra\subset\cO_X$,
the log canonical threshold of $\fra$ is denoted by
$\lct(\fra)$. 
For the precise definition in terms of a log resolution of $(X, \fra)$, we refer to
\S\ref{mult_ideals}. Given a point $x\in V(\fra)$, one defines $\lct_x(\fra)$ to be 
$\lct(\fra\vert_U)$, where $U$ is a small enough open neighborhood of $x$ in $X$.

On the other hand, suppose that $W$ is a smooth scheme of finite type over a perfect
field $L$ of positive characteristic $p$. For
a proper ideal $\fra\subset \cO_W$, the $F$-pure threshold
$\fpt(\fra)$ was introduced and studied in \cite{TW}. Given $x\in V(\fra)$, one defines as before the local version of this invariant, denoted $\fpt_x(\fra)$.
The original definition of the $F$-pure threshold involved notions and constructions from tight closure theory. However,
since we always assume that the ambient scheme is smooth, one can use an alternative
description, following \cite{MTW} and \cite{BMS2} (see \S\ref{mult_ideals} below).
 Part of the interest in the study of the $F$-pure threshold comes  from the fact that it shares many of the formal properties of the log canonical threshold. 

Before stating our main result, let us recall the
fundamental connection between log canonical thresholds and $F$-pure thresholds
via reduction mod $p$. Suppose that $X$ and $\fra\subset\cO_X$ are defined over $k$, as above.
We may choose a subring $A\subset k$, finitely generated over $\ZZ$, and models 
 $X_A$ and $\fra_A\subset\cO_{X_A}$ for $X$ and respectively $\fra$, defined over $A$. In particular, given
any closed point $s\in\Spec A$, we may consider the corresponding reductions $X_s$ and
$\fra_s\subset \cO_{X_s}$ defined over the finite residue field of $s$ denoted $k(s)$. One of the main results in \cite{HY}
implies the following relation between log canonical thresholds and $F$-pure thresholds:
after possibly replacing $A$ by a localization $A_a$, for some nonzero $a\in A$, the following hold:
\begin{enumerate}
\item[i)] $\lct(\fra)\geq\fpt(\fra_s)$ for every closed point in $s\in\Spec A$.
\item[ii)] There is a sequence of closed points $s_m\in\Spec A$ with
$\lim_{m\to\infty}{\rm char}(k(s_m))=\infty$ and such that
$\lim_{m\to\infty}\fpt(\fra_{s_m})=\lct(\fra)$.
\end{enumerate}
It is worth pointing out that a fundamental open problem in the field predicts that in this setting
there is a dense set of closed points $S\subset\Spec A$ such that
$\lct(\fra)=\fpt(\fra_s)$ for every $s\in S$. 

We now turn to the description of our main result. For every $n\geq 1$, let
${\mathcal Lct}_n$ be the set of all $\lct(\fra)$, where the pair $(X,\fra)$ is as above, with
$\dim(X)=n$. Similarly, given $n$ and a prime $p$, let ${\mathcal Fpt}(p)_n$ be the set of all
$\fpt(\fra)$, where $(W,\fra)$ is as above, with $\dim(W)=n$, and $W$ defined over a field of characteristic $p$. The following is our main result.

\begin{theorem}\label{main_intro}
For every $n\geq 1$,
the set of limit points of all sequences $(c_p)$, where $c_p\in{\mathcal Fpt}(p)_n$ 
for every prime $p$,
coincides with ${\mathcal Lct}_n$.
\end{theorem}

A key ingredient in the proof of Theorem~\ref{main_intro} is provided by 
ultraproduct constructions. Note that if $c\in {\mathcal Lct}_n$ is given as
$c=\lct(\fra)$, then the above mentioned
results in \cite{HY} (more precisely, property ii) above) imply that $c=\lim_{p\to\infty}c_p$, where 
for $p\gg 0$ prime, $c_p$ is the $F$-pure threshold of a suitable reduction $\fra_s\subset
\cO_{X_s}$,
with ${\rm char}(k(s))=p$. Thus the interesting statement in the above theorem is the converse:
given pairs $(W_m,\fra_m)$ over $L_m$ with $\dim(W_m)=n$, 
$\lim_{m\to\infty}{\rm char}(L_m)=\infty$, 
and with $\lim_{m\to\infty}\fpt(\fra_m)=c$, there is a pair $(X,\fra)$ in characteristic zero, 
with $\dim(X)=n$ and such that $c=\lct(\fra)$.

 It is easy to see that we may assume
that each $W_m=\Spec(L_m[x_1,\ldots,x_n])$,
and $c_m=\fpt_0(\fra_m)$ for some $\fra_m\subseteq (x_1,\ldots,x_n)$.
If we put $\fra_m^{(d)}=\fra_m+(x_1,\ldots,x_n)^d$, then we have $|\fpt(\fra_m^{(d)})-\fpt(\fra_m)|
\leq\frac{n}{d}$ for all $m$ and $d$. 
Ultraproduct constructions give non-standard extensions of our algebraic structures. 
In particular, we get a field $k=[L_m]$ of characteristic zero.
Since all ideals $\fra_m^{(d)}$ are generated in degree $\leq d$, they determine an ideal
$\fra^{(d)}$ in $k[x_1,\ldots,x_n]$. The key point is to show that for every $\epsilon>0$, we have
$|\lct_0(\fra^{(d)})-\fpt_0(\fra_m)|<\epsilon$ for infinitely many $m$. 
This easily implies that $\lim_{m\to\infty}\lct_0(\fra^{(d)})=c$, and since 
${\mathcal Lct}_n$ is closed by \cite[Theorem~1.3]{dFM} (incidentally, this is proved in
\emph{loc. cit.} also by non-standard arguments), we conclude that $c\in {\mathcal Lct}_n$. 

As in \cite{HY}, the result relating the log canonical threshold of $\fra^{(d)}$ and the $F$-pure thresholds
of $\fra_m^{(d)}$ follows from a more general
result relating the multiplier ideals of $\fra^{(d)}$ and the test ideals of $\fra_m^{(d)}$ (see Theorem~\ref{main_multiplier}
below). We prove this by following, with some simplifications, the main line of argument in \cite{HY}
in our non-standard setting.

The use of ultraproduct techniques in commutative algebra has been pioneered by
Schoutens (see \cite{Schoutens} and the list of references therein). 
This point of view has been particularly effective for passing
from positive characteristic to characteristic zero in an approach to tight closure theory and to its applications. Our present work combines ideas of Schoutens \cite{Sch05}
with the non-standard approach to studying limits of log canonical thresholds and
$F$-pure thresholds from 
\cite{dFM} and, respectively,
\cite{BMS1}. 

The paper is structured as follows. In \S\ref{mult_ideals} we review the definitions of multiplier ideals and test ideals, and recall how the log canonical threshold
and the $F$-pure threshold appear as the first jumping numbers in these families of ideals.
In \S\ref{non_standard} we review the basic definitions involving ultraproducts. For the 
benefit of the reader, we also
describe in detail how to go from schemes, morphisms, and sheaves
over an ultraproduct of fields to sequences
of similar objects defined over the corresponding fields. The proof of Theorem~\ref{main_intro} is given 
in \S\ref{limits}.

\subsection*{Acknowledgment}
Our project started during the "Mathematics Research Community" on Commutative Algebra,
in Snowbird, Utah, 2010. We are grateful to the AMS for setting up this program, and for providing
a stimulating research environment.

\section{Multiplier ideals and test ideals}\label{mult_ideals}

In this section we review the basic facts that we will need about multiplier ideals and test ideals.
Both these concepts can be defined under mild assumptions on the singularities of the ambient space. However, since our main result only deals with smooth varieties, we will restrict to this setting in order to simplify the definitions.

\subsection{Multiplier ideals and the log canonical threshold}\label{subsection_mult_ideals}
In what follows we recall the definition and some basic properties of multiplier ideals and log canonical thresholds. For details and further properties, we refer the reader to 
\cite[\S 9]{positivity}.

Let $k$ be a field of characteristic zero, and $X$ an irreducible and smooth scheme of finite type over $k$. Given a nonzero ideal\footnote{Every ideal sheaf that we consider is assumed to be coherent.} $\fra$  on $X$, its multiplier ideals are defined as follows.
Let us fix a log resolution of the pair $(X,\fra)$: this is a projective, birational morphism
$\pi\colon Y\to X$, with $Y$ smooth, and $\fra\cdot\cO_Y=\cO_Y(-F)$ for an effective divisor $F$,
such that $F+{\rm Exc}(\pi)$ is a divisor with simple normal crossings. Here ${\rm Exc}(\pi)$
denotes the exceptional divisor of $\pi$. 
Such resolutions exist by Hironaka's theorem, since we are in characteristic zero.
Recall that $K_{Y/X}$ denotes the relative canonical divisor
of $\pi$: this is an effective divisor supported on ${\rm Exc}(\pi)$ such that
$\cO_Y(K_{Y/X})\simeq\omega_Y\otimes f^*(\omega_X)$.
With this notation, the multiplier ideal of $\fra$ of exponent $\lambda\in\RR_{\geq 0}$ 
is defined by
\begin{equation}\label{eq1_section1}
\cJ(\fra^{\lambda}):=\pi_*\cO_Y(K_{Y/X}-\lfloor\lambda F\rfloor).
\end{equation}
Here, for a divisor with real coefficients $E=\sum_ia_iE_i$, we write
$\lfloor E\rfloor=\sum_i\lfloor a_i\rfloor E_i$, where $\lfloor a_i\rfloor$ is the largest integer $\leq a_i$.
It is a basic fact that the definition of multiplier ideals is independent of resolution. 

Let us consider some easy consequences of the definition (\ref{eq1_section1}).
If $\lambda<\mu$, then $\cJ(\fra^{\mu})\subseteq\cJ(\fra^{\lambda})$. Furthermore, given $\lambda$,
there is $\epsilon>0$ such that $\cJ(\fra^{\lambda})=\cJ(\fra^{\mu})$ whenever
$\lambda\leq\mu\leq\lambda+\epsilon$. A positive $\lambda$ is a \emph{jumping number}
of $\fra$ if $\cJ(\fra^{\lambda})\neq\cJ(\fra^{\mu})$ for all $\mu<\lambda$. If 
we write $F=\sum_ia_iE_i$, it follows from (\ref{eq1_section1})
that if $\lambda$ is a jumping number, then $\lambda a_i\in\ZZ$ for some $i$. 
In particular, we see that the jumping numbers of $\fra$ form a discrete set of rational numbers.

Suppose now that $\fra\neq\cO_X$. 
The smallest jumping number of $\fra$ is the \emph{log canonical threshold} $\lct(\fra)$. 
Note that if $0\leq\lambda\ll 1$, then $\cJ(\fra^{\lambda})=\cO_X$, hence
$\lct(\fra)=\min\{\lambda\mid\cJ(\fra^{\lambda})\neq\cO_X\}$ (this is 
finite since $\fra\neq\cO_X$). If $\fra\subseteq\frb$, then
$\cJ(\fra^{\lambda})\subseteq\cJ(\frb^{\lambda})$ for all $\lambda$; in particular,
we have $\lct(\fra)\leq\lct(\frb)$. We make the convention $\lct(0)=0$ and $\lct(\cO_X)=\infty$.

It is sometimes convenient to also have available a local version of the log canonical threshold.
If $x\in X$, then we put
$\lct_x(\fra):=\max_V\lct(\fra\vert_V)$, where the maximum ranges over all open neighborhoods
$V$ of $x$. Equivalently, we have
$$\lct_x(\fra)=\min\{\lambda\mid\cJ(\fra^{\lambda})\cdot\cO_{X,.x}\neq
\cO_{X,x}\}$$
(with the convention that this is $0$ if $\fra=(0)$, and it is infinite if $x\not\in V(\fra)$).
Note that given a proper ideal $\fra$ on $X$, there is a closed point $x\in X$ such that $\lct(\fra)=\lct_x(\fra)$.

The definition of multiplier ideals commutes with extension of the base field, as follows.
For a proof, see the proof of  \cite[Propositions~2.9]{dFM}.

\begin{proposition}\label{prop0_section1}
Let $\fra$ be an ideal on $X$. If $k\subset k'$ is a field extension, and $\phi\colon X'=X\times_
{\Spec k}\Spec k'\to
X$
and $\fra'=\fra\cdot\cO_{X'}$, then $\cJ({\fra'}^{\lambda})=\cJ(\fra^{\lambda})\cdot\cO_{X'}$
for every $\lambda\in\RR_{\geq 0}$. In particular, $\lct_{x'}(\fra')=\lct_{\phi(x')}(\fra)$ for every
$x'\in X'$.
\end{proposition}

Recall from \S\ref{intro} that ${\mathcal Lct}_n$ consists of all nonnegative rational numbers of the form $\lct(\fra)$, where $\fra$ is a proper ideal on an $n$-dimensional smooth projective variety
over a field $k$ of characteristic zero. It is clear that equivalently, we may consider the invariants
$\lct_x(\fra)$, where $(X,\fra)$ is as above, and $x\in X$ is a closed point. Furthermore, by
Proposition~\ref{prop0_section1} we may assume that $k$ is algebraically closed. One can show that in this definition we can fix the algebraically closed field $k$ and assume that 
$X=\AAA_k^n$,
and obtain the same set
(see \cite[Propositions~3.1 and 3.3]{dFM}). Furthermore, we will make use of the fact that ${\mathcal Lct}_n$
is a closed set (see \cite[Theorem~1.3]{dFM}).

\subsection{Test ideals and the $F$-pure threshold}\label{subsection_F_pure}
In this section we assume that $X$ is an irreducible, Noetherian, smooth scheme of characteristic $p>0$. 
We also assume that $X$ is $F$-finite, that is, the Frobenius morphism $F\colon X\to X$ is finite
(in fact, most of the time $X$ will be a scheme of finite type over a perfect field, in which case
this assumption is clearly satisfied). Recall that for an ideal $J$ on $X$, the $e^{\rm th}$ Frobenius
power $J^{[p^e]}$ is generated by $u^{p^e}$, where $u$ varies over the (local) generators of $J$.

Suppose that $\frb$ is an ideal on $X$. Given a positive integer $e$, one can show that there is a unique smallest ideal $J$ such that $\frb\subseteq J^{[p^e]}$. This ideal is denoted by 
$\frb^{[1/p^e]}$. Given a nonzero ideal $\fra$ and $\lambda\in\RR_{\geq 0}$, one has
$$(\fra^{\lceil \lambda p^e\rceil})^{[1/p^e]}
\subseteq(\fra^{\lceil \lambda p^{e+1}\rceil})^{[1/p^{e+1}]}$$
for all $e\geq 1$ (here $\lceil u\rceil$ denotes the smallest integer $\geq u$). By the Noetherian property, it follows that there is an ideal
$\tau(\fra^{\lambda})$ that is equal to $(\fra^{\lceil \lambda p^e\rceil})^{[1/p^e]}$ for all 
$e\gg 0$. This is the \emph{test ideal} of $\fra$ of exponent $\lambda$. 
For details and basic properties of test ideals, we refer to \cite{BMS2}. 

It is again clear that if $\lambda<\mu$, then $\tau(\fra^{\mu})\subseteq\tau(\fra^{\lambda})$. 
It takes a little argument to show that given any $\lambda$, there is $\epsilon>0$ such that
$\tau(\fra^{\lambda})=\tau(\fra^{\mu})$ whenever $\lambda\leq\mu\leq\lambda+\epsilon$
(see \cite[Proposition~2.14]{BMS2}).
We say that $\lambda>0$ is an $F$-\emph{jumping number} of $\fra$ if
$\tau(\fra^{\lambda})\neq\tau(\fra^{\mu})$ for every $\mu<\lambda$. 
It is proved in \cite[Theorem~3.1]{BMS2} that if $X$ is a scheme of finite type over an $F$-finite
field, then the $F$-jumping numbers of $\fra$ form a discrete set of rational numbers. 

The smallest $F$-jumping number of $\fra$ is the $F$-\emph{pure threshold} $\fpt(\fra)$. 
Since $\tau(\fra^{\lambda})=\cO_X$ for $0\leq\lambda\ll 1$, the $F$-pure threshold is characterized
by
$$\fpt(\fra)=\min\{\lambda\mid\tau(\fra^{\lambda})\neq\cO_X\}.$$
Note that this is finite if and only if $\fra\neq\cO_X$. We make the convention that
$\fpt(\fra)=0$ if $\fra=(0)$. 

We have a local version of the $F$-pure threshold: given
$x\in X$, we put $\fpt_x(\fra):=\max_V\fpt(\fra\vert_V)$, where the maximum is over all 
open neighborhoods $V$ of $x$. It can be also described by
$$\fpt_x(\fra)=\min\{\lambda\mid\tau(\fra^{\lambda})\cdot\cO_{X,x}\neq\cO_{X,x}\},$$
and it is finite if and only if $x\in V(\fra)$. Note that given any $\fra$, there is $x\in X$ such that
$\fpt(\fra)=\fpt_x(\fra)$.

We will make use of the following two properties of $F$-pure thresholds.
For proofs, see \cite[Proposition~2.13]{BMS2} and, respectively,  \cite[Corollary~3.4]{BMS1}.

\begin{proposition}\label{completion}
If $\fra$ is an ideal on $X$ and $S=\widehat{\cO_{X,x}}$ is the completion of the local ring of $X$ at 
a point $x\in X$, then
$\tau(\fra^{\lambda})\cdot S=\tau((\fra\cdot S)^{\lambda})$ for every $\lambda\geq 0$.
In particular, $\fpt_x(\fra)=\fpt(\fra\cdot S)$.
\end{proposition}

\begin{proposition}\label{prop2_section1}
If $\fra$ and $\frb$ are ideals on $X$, and $x\in V(\fra)\cap V(\frb)$ is such that
$\fra\cdot\cO_{X,x}+\frm^r=\frb\cdot\cO_{X,x}+\frm^r$ for some $r\geq 1$, 
where $\frm$ is the maximal ideal in $\cO_{X,x}$, then
$$|\fpt_x(\fra)-\fpt_x(\frb)|\leq\frac{\dim(\cO_{X,x})}{r}.$$
\end{proposition}

The local $F$-pure threshold admits the following alternative description, following
\cite{MTW}. If $\fra$ is an ideal on $X$ and $x\in V(\fra)$, let $\nu(e)$ denote the largest $r$
such that $\fra^r\cdot\cO_{X,x}\not\subseteq\frm^{[p^e]}$, where $\frm$ is the maximal ideal in
$\cO_{X,x}$ (we make the convention $\nu(e)=0$ if $\fra=0$). One can show that
\begin{equation}\label{F_threshold}
\fpt_x(\fra)=\lim_{e\to\infty}\frac{\nu(e)}{p^e}
\end{equation}
(see \cite[Proposition~2.29]{BMS2}). This immediately implies the assertion in the following proposition.

\begin{proposition}\label{field_extension}
Let $L\subset L'$ be a field extension of $F$-finite fields of positive characteristic.
If $\fra\subseteq L[x_1,\ldots,x_n]$ is an ideal vanishing at the origin, and
$\fra'=\fra\cdot L'[x_1,\ldots,x_n]$, then $\fpt_0(\fra)=\fpt_0(\fra')$.
\end{proposition}

Recall that we have introduced in \S\ref{intro} the set ${\mathcal Fpt}(p)_n$
consisting of all invariants of the form $\fpt(\fra)$, where $\fra$ is a proper ideal on an irreducible,
$n$-dimensional smooth scheme of finite type over $L$, with $L$ a perfect field of characteristic $p$. We can define two other related subsets of $\RR_{\geq 0}$. Let 
${\mathcal Fpt}(p)'_n$ be the set of invariants $\fpt_0(\fra)$, where $\fra\subset L[x_1,\ldots,x_n]$
is an ideal vanishing at the origin, and  $L$ is an algebraically closed field of characteristic $p$. 
We also put ${\mathcal Fpt}(p)''_n$ for the set of all $\fpt(\fra)$, where $\fra$ is a proper ideal on an
irreducible, smooth, $n$-dimensional $F$-finite scheme of characteristic $p$. We clearly have the following inclusions
\begin{equation}\label{inclusions}
{\mathcal Fpt}(p)'_n\subseteq {\mathcal Fpt}(p)_n\subseteq {\mathcal Fpt}(p)''_n.
\end{equation}

\begin{proposition}\label{dense}
${\mathcal Fpt}(p)'_n$ is dense in ${\mathcal Fpt}(p)''_n$ ${\rm (}$hence 
also in ${\mathcal Fpt}(p)_n$${\rm )}$.
\end{proposition}

This implies that in Theorem~\ref{main_intro} we may replace the sets ${\mathcal Fpt}(p)_n$
by ${\mathcal Fpt}(p)'_n$ or by ${\mathcal Fpt}(p)''_n$. 

\begin{proof}[Proof of Proposition~\ref{dense}]
Suppose that $\fra$ is a proper ideal on $X$, where $X$ is irreducible, smooth, $F$-finite,
$n$-dimensional, and
of characteristic $p$. Let $c=\fpt(\fra)$. We can find $x\in X$ such that
$c=\fpt_x(\fra)$. By Proposition~\ref{completion}, we have $c=\fpt(\fra\cdot \widehat{\cO_{X,x}})$.
Note that by Cohen's theorem, we have an isomorphism $\widehat{\cO_{X,x}}\simeq
L\llbracket x_1,\ldots,x_d\rrbracket$, with $L$ an $F$-finite field, and $d\leq n$. 
If $\frm$ is the maximal ideal in $\widehat{\cO_{X,x}}$ and $c_i=\fpt(\fra\cdot\widehat{\cO_{X,x}}+\frm^i)$,
then Proposition~\ref{prop2_section1} gives $c=\lim_{i\to\infty} c_i$. On the other hand, 
there are ideals
$\frb_i\subset L[x_1,\ldots,x_d]$ such that $\frb_i\cdot \widehat{\cO_{X,x}}=
\fra\cdot\widehat{\cO_{X,x}}+\frm^i$, and
another application of Proposition~\ref{completion} gives $c_i=\fpt_0(\frb_i)$.
It is easy to see that 
$c_i=\fpt_0(\frb_i\cdot L[x_1,\ldots,x_n])$ (for example, this is a consequence of 
formula (\ref{F_threshold})). 
 It now follows from Proposition~\ref{field_extension}
that $c_i=\fpt_0(\frb_i\cdot\overline{L}[x_1,\ldots,x_n])$, where $\overline{L}$ is an algebraic closure
of $L$. Therefore all $c_i$ lie in ${\mathcal Fpt}(p)'_n$, which proves the proposition.
\end{proof}

In \S\ref{limits} we will use a slightly different description of the test ideals that we now present.
More precisely, we give a different description of $\frb^{[1/p^e]}$, when $\frb$ 
is an arbitrary ideal on $X$.
Suppose that $X$ is an irreducible, smooth scheme of finite type over a perfect field $L$ of characteristic $p$. Let $\omega_X=\wedge^n\Omega_{X/L}$, where $n=\dim(X)$.
Recall that the Cartier isomorphism (see \cite{DI}) gives in particular
an isomorphism $\omega_X\simeq {\mathcal H}^n(F_*\Omega_{X/L}^{\bullet})$, where 
$F$ is the (absolute) Frobenius morphism, and $\Omega^{\bullet}_{X/L}$ is the de Rham complex
of $X$. In particular, we get a surjective $\cO_X$-linear map
$t_X\colon F_*\omega_X\to\omega_X$. This can be explicitly described in coordinates, as follows.
Suppose that $u_1,\ldots,u_n\in\cO_{X,x}$ form a regular system of parameters, where $x\in X$
is a closed point.
We may assume that $u_1,\ldots,u_n$ are defined in an affine open neighborhood $U$ of $x$,
and that $du=du_1\wedge\ldots\wedge du_n$ gives a basis of $\omega_X$ on $U$.
Furthermore, we may assume that $\cO_U$ is free over $\cO_U^p$, with basis
$$\{u_1^{i_1}\cdots u_n^{i_n}\mid 0\leq i_j\leq p-1\,\text{for}\, 1\leq j\leq n\}$$
(note that the residue field of $\cO_{X,x}$ is a finite extension of $L$, hence it is perfect).
In this case $t_X$ is characterized by the fact that $t_X(h^p w)=h\cdot t_X(w)$ for every
$h\in\cO_X(U)$, and on the above basis over $\cO_X(U)^p$ it is described by
\begin{equation}\label{formula_t}
t_X(u_1^{i_1}\cdots u_n^{i_n}du)=\left\{
\begin{array}{cl}
du, & \text{if}\,\,i_j=p-1\,\text{for all}\, \,j; \\[2mm]
0, & \text{otherwise}. \\[2mm]
\end{array}\right.
\end{equation}

Iterating $e$ times $t_X$  gives $t_X^e\colon F^e_*\omega_X\to\omega_X$. These maps are functorial
in the following sense. If $\pi\colon Y\to X$ is a proper birational morphism between
irreducible smooth varieties as above, then we have a commutative diagram
\begin{equation}\label{diag3_1}
\begin{CD}
\pi^*(F^e_*(\omega_X))@>{\pi^*(t_{X}^e)}>>\pi^*(\omega_X)\\
@VVV @VV{\psi}V\\
F^e_*(\omega_Y)@>{t_Y^e}>>\omega_Y,
\end{CD}
\end{equation}
where  $\psi$ is the canonical morphism induced by pulling-back $n$-forms,
and
the left vertical map is the composition
$$\pi^*(F^e_*(\omega_X))\to F^e_*(\pi^*(\omega_X))\overset{F^e_*(\psi)}\to
F^e_*(\omega_Y).$$

Suppose now that $X$ is as above, and $\frb$ is an ideal on $X$. Since $\omega_X$
is a line bundle, it follows that the image of $F_*(\frb\cdot\omega_X)$ by $t_X^e$
can be written as $J\cdot \omega_X$, for a unique ideal $J$ on $X$. It is an easy consequence of 
the description of $\frb^{[1/p^e]}$ in \cite[Proposition~2.5]{BMS2} and of formula (\ref{formula_t})
that in fact $J=\frb^{[1/p^e]}$ (see also \cite[Proposition~3.10]{BSTZ}).

\section{A review of non-standard constructions}\label{non_standard}

We begin by reviewing some general facts about ultraproducts. For a detailed introduction to this topic,
the reader is referred to \cite{Gol}. We then explain how geometric objects over an ultraproduct
of fields correspond to sequences of such geometric objects over the fields we are starting with,
up to a suitable equivalence relation. Most of this material 
is well-known to the experts, and can be found, for example, in
\cite[\S 2]{Sch05}. However, we prefer to give a detailed presentation for the benefit of those readers having little or no familiarity with non-standard constructions.

\subsection{Ultrafilters and ultraproducts}
Recall that an \emph{ultrafilter} on the set of positive integers $\NN$ is a nonempty collection $\cU$
of subsets of $\NN$ that satisfies the following properties:
\begin{enumerate}
\item[(i)] If $A$ and $B$ lie in $\cU$, then $A\cap B$ lies in $\cU$.
\item[(ii)] If $A\subseteq B$ and $A$ is in $\cU$, then $B$ is in $\cU$.
\item[(iii)] The empty set does not belong to $\cU$.
\item[(iv)] Given any $A\subseteq\NN$, either $A$ or $\NN\smallsetminus A$ lies in $\cU$.
\end{enumerate}
An ultrafilter $\cU$ is \emph{non-principal} if no finite subsets of $\NN$ lie in $\cU$. It is an easy consequence of Zorn's Lemma that non-principal ultrafilters exist, and we fix one such ultrafilter 
$\cU$. Given a property ${\mathcal P}(m)$, where $m\in\NN$, we say that ${\mathcal P}(m)$
holds \emph{for almost all} $m$ if $\{m\in\NN\mid {\mathcal P}(m)\,\text{holds}\}$ lies in $\cU$. 

Given a sequence of sets $(A_m)_{m\in\NN}$, the \emph{ultraproduct}
$[A_m]$ is the quotient of $\prod_{m\in\NN}A_m$ by the equivalence relation 
given by  $(a_m)\sim(b_m)$ if $a_m=b_m$ for almost all $m$. 
We write the class of $(a_m)$ in $[A_m]$ by $[a_m]$. 
Note that the element $[a_m]$ is well-defined even if $a_m$ is only defined for almost all $m$. 
Similarly, the set $[A_m]$ is well-defined if we give $A_m$ for almost all $m$.

If $A_m=A$ for all $m$,
then one writes $\ltightstar{A}$ instead of $[A_m]$. This is the \emph{non-standard extension}
of $A$. Note that there is an obvious inclusion $A\hookrightarrow\ltightstar{A}$, that takes $a\in A$
to the class of the constant sequence $(a)$.

The general principle is that
if all $A_m$ have a certain algebraic structure, then so does $[A_m]$, by defining the corresponding structure
component-wise on $\prod_{m\in\NN}A_m$. For example, if we consider fields $(L_m)_{m\in\NN}$, then 
$k:=[L_m]$ is a field. 
In particular, the non-standard extension $\ltightstar{\RR}$ of $\RR$ is an ordered field. 
Furthermore,
it is easy to see that if all $L_m$ are algebraically closed, then so is $k$. Note also that if
$\lim_{m\to\infty}{\rm char}(L_m)=\infty$, then ${\rm char}(k)=0$. 

Given a sequence of maps $f_m\colon A_m\to B_m$, for $m\in\NN$, we get a map
$[f_m]\colon [A_m]\to [B_m]$ that takes $[a_m]$ to $[f_m(a_m)]$. 
In particular, given a map $f\colon A\to B$, we get a map $\ltightstar{f}
\colon\ltightstar{A}\to\ltightstar{B}$ that extends $f$.
 If each $A_m$ is a subset of $B_m$, we can identify $[A_m]$ to a subset of $[B_m]$
 via the corresponding map. The subsets of $[B_m]$ of this form are called \emph{internal}.

We will use in \S\ref{limits} the following notion. Suppose that  $u=[u_m]\in\ltightstar{\RR}$ 
is bounded
(this means that there is $M\in\RR_{>0}$ such that $\ltightstar{|u|}\leq M$, that is,
$|u_m|\leq M$ for almost all $m$). In this case, there is a unique
real number, the \emph{shadow} ${\rm sh}(u)$ of $u$, with the property that for every positive 
real number $\epsilon$, we have $\ltightstar{|u-{\rm sh}(u)|}<\epsilon$, that is,
$|{\rm sh}(u)-u_m|<\epsilon$ for almost all $m$. We refer to 
\cite[\S 5.6]{Gol} for a discussion of shadows. A useful property is that if $(c_m)_{m\in\NN}$
is a convergent sequence, with $\lim_{m\to\infty}c_m=c$, then ${\rm sh}([c_m])=c$
(see \cite[Theorem~6.1]{Gol}). On the other hand, it is a consequence of the definition that ${\rm sh}([c_m])$ is the limit of a suitable subsequence of  $(c_m)_{m\in\NN}$.

\subsection{Schemes, morphisms, and sheaves over an ultraproduct of fields}
Suppose that $\cU$ is a non-prinicpal ultrafilter on $\NN$ as in the previous section, and
suppose that $(L_m)_{m\in\NN}$ is a sequence of fields. We denote the corresponding
ultraproduct by $k=[L_m]$. Let us temporarily fix $n\geq 1$, and consider the polynomial rings
$R_m=L_m[x_1,\ldots,x_n]$. We write 
$k[x_1,\ldots,x_n]_{\rm int}$ for the ring $[R_m]$,
the ring of \emph{internal polynomials} in $n$ variables (we emphasize, however, that the elements of this ring are
\emph{not} polynomials). Given a sequence of ideals $(\fra_m\subseteq R_m)_{m\in\NN}$,
we get the \emph{internal ideal} $[\fra_m]$ in $k[x_1,\ldots,x_n]_{\rm int}$.

We have an embedding $k[x_1,\ldots,x_n]\hookrightarrow
k[x_1,\ldots,x_n]_{\rm int}$.
Its image consists of those
$g=[g_m]\in k[x_1,\ldots,x_n]_{\rm int}$ for which there is $d\in\NN$ such that $\deg(g_m)\leq d$ for almost all $m$
(in this case we say that $g$ \emph{has bounded degree}). We say that an ideal
$\frb\subseteq k[x_1,\ldots,x_n]_{\rm int}$ is generated in bounded degree if it is generated by an ideal in $k[x_1,\ldots,x_n]$
(in which case $\frb$ is automatically an internal ideal). Given an ideal $\fra$ in $k[x_1,\ldots,x_n]$,
we put $\fra_{\rm int}:=\fra\cdot k[x_1,\ldots,x_n]_{\rm int}$.

The connection between $k[x_1,\ldots,x_n]$ and $k[x_1,\ldots,x_n]_{\rm int}$ is studied in 
\cite{vdDS}. In particular, the following is \cite[Theorem~1.1]{vdDS}.

\begin{theorem}\label{thm1_vdDS}
The extension $k[x_1,\ldots,x_n]\hookrightarrow k[x_1,\ldots,x_n]_{\rm int}$ is faithfully flat.
In particular, given any ideal $\fra$ in $k[x_1,\ldots,x_n]$, we have $\fra_{\rm int}
\cap k[x_1,\ldots,x_n]=\fra$.
\end{theorem}

It follows from the theorem that ideals of $k[x_1,\ldots,x_n]_{\rm int}$ generated in bounded degree are in order-preserving bijection with the ideals in $k[x_1,\ldots,x_n]$. Furthermore, note that every such ideal of
$k[x_1,\ldots,x_n]_{\rm int}$ is of the form $[\fra_m]$ for a sequence $(\fra_m)_{m\in\NN}$ that is \emph{generated in bounded degree}, that is, such that for some $d$,  
 $\fra_m\subseteq L_m[x_1,\ldots,x_n]$ is generated by polynomials of degree $\leq d$ for almost all $m$. 
 Of course, we have $[\fra_m]=[\frb_m]$ if and only $\fra_m=\frb_m$ for almost all $m$.
 Given such a sequence $(\fra_m)_{m\in\NN}$, we call $[\fra_m]\cap k[x_1,\ldots,x_n]$ the
 \emph{ideal of polynomials} corresponding to the sequence.
 
 Our next goal is to describe how to associate to a geometric object over $k$ a sequence of
 corresponding objects over each of $L_m$ (in fact, an equivalence class of such sequences). 
 Given a separated scheme $X$ of finite type over $k$, we will associate to it an
 \emph{internal scheme} $[X_m]$, 
 by which we mean the following: we have schemes $X_m$ of finite type over $L_m$ for almost all $m$; furthermore, two such symbols $[X_m]$ and $[Y_m]$ define the same equivalence class
 if $X_m=Y_m$ for almost all $m$. An  \emph{internal morphism}
 $[f_m]\colon [X_m]\to [Y_m]$ between internal schemes consists of an equivalence class
 of sequences of morphisms of schemes $f_m\colon X_m\to Y_m$ (defined for almost all $m$), where $[f_m]=[g_m]$
 if $f_m=g_m$ for almost all $m$.
 
 We want to define a functor $X\to X_{\rm int}$ from separated schemes of finite type over $k$ to 
 internal schemes. We first consider the case when $X$ is affine.
 In this case let us choose a closed embedding
 $X\hookrightarrow \AAA_k^N$, defined by the ideal $\fra\subseteq k[x_1,\ldots,x_N]$. If
 $\fra_{\rm int}=[\fra_m]$, then we take $X_m$ to be defined in
 $\AAA_{L_m}^N$ by $\fra_m$. 
 Note that $X_{\rm int}:=[X_m]$ is well-defined. 
 We also put $\cO(X)_{\rm int}:=[L_m[x_1,\ldots,x_N]/\fra_m]$, and note that we have a canonical ring
 homomorphism $\eta_X\colon \cO(X)\to \cO(X)_{\rm int}$.
 Suppose now that we have a morphism $f\colon Y\to X$ of affine schemes as above, and closed embeddings $Y\hookrightarrow \AAA_k^N$ and $X\hookrightarrow\AAA_k^M$. 
 We have a homomorphism $\phi\colon k[x_1,\ldots,x_M]\to k[x_1,\ldots,x_N]$ that induces $f$,
 and that extends to an internal morphism $k[x_1,\ldots,x_M]_{\rm int}
 \to k[x_1,\ldots,x_N]_{\rm int}$. This induces morphisms $f_m\colon Y_m\to X_m$ 
 for almost all $m$, hence an internal morphism $Y_{\rm int}\to X_{\rm int}$. 
 It is easy to see that this is independent of the choice of the lifting $\phi$, and that it is functorial.
 
 The first consequence is that if we replace $X\hookrightarrow \AAA_k^N$ by a different embedding
 $X\hookrightarrow\AAA_k^M$, then the two internal schemes that we obtain are canonically isomorphic. 
 We use this to extend the above definition to the case when $X$ is not necessarily affine,
 as follows. Note first that if $\overline{L_m}$ is an algebraic closure of $L_m$, and if $\overline{k}=
 [\overline{L_m}]$, then $\overline{k}$ is an algebraically closed field containing $k$, and for every
 affine $X$ as above, with $X_{\rm int}=[X_m]$, we have a natural bijection of sets
 $X(\overline{k})\simeq [X_m(\overline{L_m})]$. 
 
 \begin{lemma}\label{open_immersion}
 Let $X$ be an affine scheme as above, $U\subset X$ an affine open subset, and write
 $X_{\rm int}=[X_m]$ and $U_{\rm int}=[U_m]$.
 \item[(i)] The induced maps $U_m\to X_m$ are open immersions for almost all $m$.
 \item[(ii)] If $X=U_1\cup\ldots\cup U_r$ is an affine open cover, and
 $(U_i)_{\rm int}=[U_{i,m}]$ for every $i$, then $X_m=U_{1,m}\cup\ldots\cup U_{r,m}$
 for almost all $m$.
 \end{lemma} 
 
 \begin{proof}
 The first assertion is clear in the case when $U$ is a principal affine open subset 
 corresponding to $f\in \cO(X)$: if the image of $f$ in $\cO(X)_{\rm int}$ is $[f_m]$, then for almost all $m$ we have that $U_m$ is the principal affine open subset of $X_m$ corresponding to $U_m$. 
 The assertion in (ii) is clear, too, when all $U_i$ are principal affine open subsets in $X$:
 once we know that the $U_{i,m}$ are open in $X_m$, to get the assertion we want it is enough to 
look at the $\overline{k}$-valued points of $X$.

We now obtain the assertion in (i) in general, since we may cover $U$ by finitely many principal
affine open subsets in $X$ (hence also in $U$). We then deduce (ii) in general from (i) by considering the $\overline{k}$-valued points of $X$.
\end{proof}

 Given any scheme $X$, separated and of finite type over $k$, consider an affine open cover 
 $X=U_1\cup\ldots\cup U_r$, and let $(U_i)_{\rm int}=[U_{i,m}]$. The intersection
 $U_i\cap U_j$ is affine and open in both $U_i$ and $U_j$, hence by Lemma~\ref{open_immersion},
 $U_{i,m}\cap U_{j,m}$ is affine and open in both $U_{i,m}$ and $U_{j,m}$ for almost all $m$.
 We get $X_m$ by glueing, for all $i$ and $j$, the open subsets $U_{i,m}$ and $U_{j,m}$ along $(U_i\cap U_j)_m$, and put
 $X_{\rm int}=[X_m]$. It is straightforward to check that $X_{\rm int}$ is independent of the choice of cover (up to a canonical isomorphism). 
 Similarly, given a morphism 
 of schemes $f\colon Y\to X$ we get an internal morphism $f_{\rm int}=[f_m]\colon 
 Y_{\rm int}\to X_{\rm int}$
 by gluing the internal morphisms obtained by restricting $f$ to suitable affine open subsets.
 Therefore we have a functor from the category of separated schemes of finite type over $k$ to the category
 of internal schemes and internal morphisms. This has the property that given $L_m$-algebras 
 $A_m$ for almost all $m$, if
 $A=[A_m]$, then we have a natural bijection of sets
 \begin{equation}\label{eq_points}
 {\rm Hom}(\Spec A, X)\simeq [{\rm Hom}(\Spec A_m, X_m)],
 \end{equation} 
 where $X_{\rm int}=[X_m]$.
 In particular, we have a bijection $X(\overline{k})\simeq [X_m(\overline{L_m})]$.
 
 We do not attempt to give a comprehensive account of the properties of this construction,
 but list in the following proposition a few that we will need.

\begin{proposition}\label{general_properties}
Let $X$ and $Y$ be separated schemes of finite type over $k$, and $X_{\rm int}=[X_m]$ and 
$Y_{\rm int}=[Y_m]$
the corresponding internal schemes.
\item[(i)] For every affine open subset $U$ of $X$, the ring homomorphism
$\eta_U\colon \cO_X(U)\to\cO_X(U)_{\rm int}$ is faithfully flat.
\item[(ii)] $X$ is reduced or integral if and only if  $X_m$ has the same property for almost all
 $m$.
\item[(iii)] The internal scheme corresponding to $X\times Y$ is $[X_m\times Y_m]$. 
\item[(iv)] If $f\colon Y\to X$ is an open or closed immersion, then the induced morphisms
$f_m\colon Y_m\to X_m$ have the same property for almost all $m$. In particular, $X_m$
is separated for almost all $m$.
\item[(v)] If $X^{(1)},\ldots,X^{(r)}$ are the irreducible components of $X$, 
and $X^{(i)}_{\rm int}=[X^{(i)}_m]$, then $X^{(1)}_m,\ldots,X^{(r)}_m$ are the irreducible components
of $X_m$ for almost all $m$.
\item[(vi)] If $X$ is affine and $f\colon Y\to X$ is a projective morphism, then $f_m\colon Y_m
\to X_m$ is projective for almost all $m$.
\item[(vii)] If $\dim(X)=d$, then $\dim(X_m)=d$ for almost all $m$.
\end{proposition}

\begin{proof}
The assertion in (i) follows from definition and Theorem~\ref{thm1_vdDS}. The assertions
in (ii) follow from definition and the fact that if $\fra$ is an ideal in $k[x_1,\ldots,x_N]$, then
$\fra$ is prime or radical if and only if $\fra\cdot k[x_1,\ldots,x_n]_{\rm int}$ has the same property,
see \cite[Theorem~2.5, Corollary~2.7]{vdDS}. Properties (iii) and (iv) are easy consequences of the definition (note that we have already checked the assertion regarding open immersions when both $X$ and $Y$ are affine). The second assertion in (iv) follows from the fact that the diagonal map 
$X\to X\times X$ being a closed immersion implies that $X_m\to X_m\times X_m$ is a closed immersion for almost all $m$. We obtain (v) from (ii), (iv), and the fact that
$X_m^{(1)},\ldots,X_m^{(r)}$ cover $X_m$ for almost all $m$. This follows by computing
the $\overline{k}$-points
of $X$, and using (\ref{eq_points})

 In order to prove (vi), note that 
$(\PP_k^N)_{\rm int}\simeq [\PP^N_{L_m}]$. Therefore a closed embedding $\iota\colon Y\hookrightarrow X\times
\PP_k^N$ induces by (iii) and (iv) closed embeddings $\iota_m\colon Y_m\hookrightarrow
X_m\times\PP_{L_m}^N$ for almost all $m$. 

We prove (vii) by induction on $\dim(X)$. Using v), we reduce to the case when $X$ is irreducible.
After replacing $X$ by $X_{\rm red}$, we see that we may assume, in fact, that
$X$ is integral, hence by (ii), for almost all $m$ we have $X_m$ integral. It is enough to prove the assertion for an affine open subset $U$ of $X$, hence we may assume that $X=\Spec A$ is affine,
and let us write $X_m=\Spec A_m$.
If $f\in\cO(X)$ is nonzero, then $\dim(A/(f))=\dim(A)-1$. Let $[f_m]=\eta_X(f)\in
\cO(X)_{\rm int}$, hence for almost all $m$ we have $f_m\neq 0$ and
$\dim(A_m/(f_m))=\dim(A_m)-1$. Since the internal scheme corresponding to $\Spec A/(f)$
is $[\Spec A_m/(f_m)]$, we conclude by induction.
\end{proof}

Suppose now that $X$ is a scheme over $k$ as above, and $\cF$ is a coherent sheaf on $X$.
If $X_{\rm int}=[X_m]$, we define an internal coherent sheaf on
$[X_m]$ to be a symbol $[\cF_m]$, where $\cF_m$ is defined for almost all $m$, and it is a coherent sheaf of $X_m$. Furthermore, two such symbols $[\cF_m]$ and $[\cF'_m]$ are identified precisely 
when $\cF_m=\cF'_m$ for almost all $m$. A morphism of internal coherent sheaves is defined 
in a similar way, and we get an abelian category consisting of internal coherent sheaves on 
$[X_m]$. 

We now define a functor $\cF\to\cF_{\rm int}$  from the category of coherent sheaves on $X$ to that of internal
coherent sheaves on $X_{\rm int}$. Given an affine open subset $U$ of $X$ and the corresponding internal scheme $U_{\rm int}=[U_m]$, we consider
the $\cO_X(U)_{\rm int}$-module $\cT_U:=\cF(U)\otimes_{\cO_X(U)}\cO_X(U)_{\rm int}$.
We claim that this is equal to $[M_m]$ for suitable $\cO_{X_m}(U_m)$-modules $M_m$. 
Indeed, this follows by considering a finite free presentation 
$$\cO_X(U)^{\oplus r}\overset{\phi}\to\cO_X(U)^{\oplus s}\to\cF(U)\to 0.$$
If $\phi$ is defined by a matrix $(a_{i,j})_{i,j}$ and if we write $\eta_U(a_{i,j})=[a_{i,j,m}]$, then 
we may take each $M_m$ to be the cokernel of the map $\cO_{X_m}(U_m)^{\oplus r}
\to\cO_{X_m}(U_m)^{\oplus s}$ defined by the matrix $(a_{i,j,m})_{i,j}$.
We put $\cF_m(U)=M_m$ for almost all $m$. It is now easy to see that the $\cF_m(U)$ glue
together for almost all $m$ to give coherent sheaves $\cF_m$ on $X_m$. 
Therefore we get an internal coherent sheaf $\cF_{\rm int}$ on $X_{\rm int}$.
Given a morphism of 
coherent sheaves on $X$, we clearly get a corresponding morphism of internal coherent sheaves. 
It follows from definition and Proposition~\ref{general_properties} (i) that this 
functor is exact in a strong sense:
a bounded complex of coherent sheaves on $X$ is acyclic if and only if 
the corresponding complexes of coherent sheaves on $X_m$ are acyclic for almost
all $m$. Note also that the functor is compatible with tensor product: if $\cF_{\rm int}=[\cF_m]$
and $\cG_{\rm int}=[\cG_m]$, then 
$(\cF\otimes_{\cO_X}\cG)_{\rm int}$ is canonically isomorphic to $[\cF_m\otimes_{\cO_{X_m}}
\cG_m]$.
We collect in the following proposition a few other properties of this functor that we will need.

\begin{proposition}\label{various_properties2}
Let $X$ be a separated scheme of finite type over $k$, and $\cF$ a coherent sheaf on $X$.
Consider $X_{\rm int}=[X_m]$ and $\cF_{\rm int}=[\cF_m]$.
\begin{enumerate}
\item[(i)] $\cF$ is locally free of rank $r$ if and only if $\cF_m$ has the same property for almost all $m$.
\item[(ii)] If $\cF$ is an ideal in $\cO_X$ defining the closed subscheme $Z$ of $X$,
and $Z_{\rm int}=[Z_m]$, then $\cF_m$ is ${\rm (}$isomorphic to${\rm )}$ the ideal defining $Z_m$ in $X_m$ for almost all $m$.
\item[(iii)] If $f\colon Y\to X$ is a morphism of schemes as above, and 
$f_{\rm int}=[f_m]$, then we have a canonical isomorphism $f^*(\cF)_{\rm int}\simeq [f_m^*(\cF_m)]$.
\item[(iv)] If $g\colon Y\to X$ is a projective morphism of schemes as above, and $g_{\rm int}=[g_m]
\colon [Y_m]\to[X_m]$, then for every $i\geq 0$ we have a canonical isomorphism 
$$R^if_*(\cF)_{\rm int}\simeq [R^i(f_m)_*(\cF_m)].$$
\item[v)] If $f$ is as in ${\rm (iv)}$, $X$ is affine, and $\cF$ is a line bundle on $X$ that is 
${\rm (}$very${\rm )}$ ample over $X$,
then $\cF_m$ is ${\rm (}$very${\rm )}$ ample over $X_m$ for almost all $m$.
\end{enumerate} 
\end{proposition}

\begin{proof}
The first assertion follows from Proposition~\ref{general_properties} (i) and the fact that given
a faithfully flat ring homomorphism $A\to B$, a finitely generated $A$-module $M$ is locally free
of rank $r$ if and only if the $B$-module $M\otimes_AB$ is locally free of rank $r$. Assertion 
(ii) is an immediate consequence of the definitions. 
In order to prove (iii) it is enough to consider the case when both $X$ and $Y$ are affine.
In this case the assertion follows from the natural isomorphism
$[M_m]\otimes_{[A_m]}[B_m]\simeq [M_m\otimes_{A_m}B_m]$ whenever 
$A_m\to B_m$ are ring homomorphisms, and the $M_m$ are finitely generated $A_m$-modules.

Let us now prove (iv). Suppose first that $X$ is affine. The first step is to  construct canonical 
morphisms
\begin{equation}\label{eq_various}
H^i(Y,\cF)_{\rm int}\to [H^i(Y_m,\cF_m)].
\end{equation}
This can be done by computing the cohomology as Cech cohomology with respect to a finite affine open cover of $Y$, and the corresponding affine open covers of $Y_m$ (and by checking that 
the definition is independent of the cover). It is enough to prove that the maps 
(\ref{eq_various}) are isomorphisms: if $X$ is not affine,
then we simply glue the corresponding isomorphisms over a suitable affine open cover of $X$. 
Since $Y$ is isomorphic to a closed subscheme of some $X\times\PP_k^N$, it is enough to prove
that the morphisms (\ref{eq_various}) are isomorphisms when $Y=\PP_X^N$. 
Explicit computation of cohomology implies that (\ref{eq_various}) is an isomorphism when
$\cF=\cO_{\PP_X^N}(\ell)$ (note that $\cO_{\PP_X^N}(\ell)_{\rm int}\simeq [\cO_{\PP^N_{X_m}}(\ell)]$).

We now prove that (\ref{eq_various}) is an isomorphism by descending induction on $i$, the
case
$i>N$ being trivial. 
Given any $\cF$, there is an exact sequence
$$0\to\cG\to\cO_{\PP_X^N}(\ell)^{\oplus r}\to\cF\to 0,$$
for some $\ell$ and $r$. We use the induction hypothesis, the long exact sequence in cohomology and
the 5-lemma to show first that (\ref{eq_various}) is surjective for all $\cF$. Applying this for $\cG$,
we then conclude that (\ref{eq_various}) is also injective for all $\cF$. This completes the proof of
(iv). The assertion in (v) follows using (iii) and Proposition~\ref{general_properties} (iv), from the fact that if $Y=\PP_X^N$ and $\cF=\cO_Y(1)$, then $Y_m\simeq \PP_{X_m}^N$ and
$\cF_m\simeq \cO_{Y_m}(1)$ for almost all $m$.
\end{proof}

We will need the following uniform version of asymptotic Serre vanishing
(see also \cite[Corollary~2.16]{Sch05}).

\begin{corollary}\label{Serre_vanishing}
Let $f\colon Y\to X$ be a projective morphism of schemes over $k$ as above, with $X$ affine.
If $\cF$ is a coherent sheaf on $Y$ and $\cL$ is a line bundle on $Y$ that is ample over $X$ and such that
$H^i(Y,\cF\otimes\cL^j)=0$ for all $i\geq 1$ and all $j\geq j_0$, then for almost all $m$ we have
$H^i(Y_m,\cF_m\otimes\cL_m^j)=0$ for all $i\geq 1$ and $j\geq j_0$, where 
$Y_{\rm int}=[Y_m]$, $\cF_{\rm int}=[\cF_m]$, and $\cL_{\rm int}=[\cL_m]$.
\end{corollary}

\begin{proof}
Note first that we may assume that $\cL$ is very ample. Indeed, if
$N$ is such that $\cL^N$ is very ample, then we may apply the very ample case to
the line bundle $\cL^N$ and to the sheaves $\cF,\cF\otimes\cL,\ldots,\cF\otimes\cL^{N-1}$ to obtain the assertion in the corollary. Let $r=\dim(Y)$. It follows from Proposition~\ref{various_properties2} (iv) that if $i\geq 1$ and
$j\geq j_0$ are fixed, then $H^i(Y_m,\cF_m\otimes\cL_m^j)=0$ for almost all $m$. 
In particular,  for almost all $m$ we have $H^i(Y_m,\cF_m\otimes \cL_m^{j})=0$ for 
$1\leq i\leq\dim(Y_m)=r$ and $j_0\leq j\leq j_0+r-1$. For every such $m$, it follows that
$\cF_m$ is $(j_0+r)$-regular in the sense of Castelnuovo-Mumford regularity, hence
$H^i(Y_m,\cF_m\otimes\cL_m^j)=0$ for every $i\geq 1$ and $j\geq j_0+r-i$
(see \cite[Chapter~1.8.A]{positivity}). This completes the proof of the corollary.
\end{proof}

\begin{proposition}\label{nonsingular}
If $X$ is a separated scheme of finite type over $k$ and $X_{\rm int}=[X_m]$, then
there is a canonical isomorphism
$(\Omega_{X/k})_{\rm int}\simeq [\Omega_{X_m/L_m}]$. In particular, 
$X$ is smooth of pure dimension $n$ if and only if $X_m$ is smooth of pure dimension $n$
for almost all $m$.
\end{proposition}

\begin{proof}
It is enough to give a canonical isomorphism $(\Omega_{X/k})_{\rm int}=[\Omega_{X_m/L_m}]$
when $X$ is affine. Note that we have such an isomorphism when $X=\AAA_k^N$. In general,
if $X$ is a closed subscheme of $\AAA_k^N$ defined by the ideal $\fra$, the sheaf
$\Omega_{X/k}$ is the cokernel of a morphism $\fra/\fra^2\to \Omega_{\AAA^N_k}\vert_X$.
If $\fra=[\fra_m]$, then
for almost all $m$ we have an analogous description of each $\Omega_{X_m/L_m}$ in terms of the embedding
$X_m\hookrightarrow\AAA_{L_m}^N$ given by $\fra_m$. Therefore we obtain the desired isomorphism, and
one can then check that this is independent of the embedding.

Recall that
$X$ is smooth of pure dimension $n$ if and only if $\dim(X)=n$ and $\Omega_{X/k}$
is locally free of rank $n$. The second assertion in the proposition now follows from the first one, together with Proposition~\ref{general_properties} (vii) and Proposition~\ref{various_properties2} (i).
\end{proof}

Suppose now that $X$ is a smooth scheme over $k$ as above, and 
$D=a_1D^{(1)}+\ldots+a_rD^{(r)}$ is a divisor on $X$. It follows from 
Proposition~\ref{general_properties} (ii), (vii) that if $D^{(i)}_{\rm int}=[D^{(i)}_{m}]$, then  
$D^{(i)}_{m}$
is a prime divisor on $X_m$ for almost all $m$. For all such $m$ we put
$D_m=a_1D^{(1)}_{m}+\ldots+a_rD^{(r)}_{m}$. 

\begin{remark}
Note that in the case when $D$ is effective, and thus can be considered as a subscheme
of $X$, the above convention is compatible with our previous definition via $D_{\rm int}=
[D_m]$. Indeed, if we define the $D_m$ via the latter formula, then it follows from definition
that since $D$ is locally defined by one nonzero element, the same holds for $D_m$ for almost all $m$. Furthermore, Proposition~\ref{general_properties} (v) implies that
$D^{(1)}_{m},\ldots,D^{(r)}_{m}$ are the irreducible components of $D_m$
for almost all $m$.
We also see that
the coefficient of $D^{(i)}_{m}$ in $D_m$ is equal to $a_i$ for almost all $m$: this follows from the fact
that this coefficient is the largest nonnegative integer $d_i$ such that $d_iD^{(i)}_{m}$
is a subscheme of $D_m$. 
\end{remark}

We thus see that for every divisor $D$, we have $\cO(D)_{\rm int}=
[\cO(D_m)]$. Indeed, when $-D$ is effective, this follows from the above remark and Proposition~\ref{various_properties2} (ii).  The general case follows easily by reducing to the case when $X$ is affine, and replacing $D$ by $D+{\rm div}(f)$ for a suitable $f\in\cO(X)$
such that $-D-{\rm div}(f)$ is effective.

\begin{proposition}\label{SNC}
Let $X$ be a smooth, separated scheme over $k$, and $D=\sum_{i=1}^ND^{(i)}$
an effective divisor on $X$, with simple normal crossings, where the $D^{(i)}$ are distinct prime divisors. If $X_{\rm int}=[X_m]$ and $D_{\rm int}=[D_m]$,
then $D_m$ has simple normal crossings for almost all $m$.
\end{proposition}

\begin{proof}
Note that $X_m$ is smooth over $L_m$ for almost all $m$ by Proposition~\ref{nonsingular}.
Since $D$ has simple normal crossings,
for every $r$ and every $1\leq i_1<\ldots<i_r\leq N$ the subscheme
$D^{(i_1)}\cap\ldots\cap D^{(i_r)}$ is smooth over $k$ (possibly empty). It follows from definition that we have
$$(D^{(i_1)}\cap\ldots\cap D^{(i_r)})_{\rm int}=[D_m^{(i_1)}\cap\ldots \cap D^{(i_r)}_m],$$
hence $D^{(i_1)}_m\cap\ldots \cap D^{(i_r)}_m$ is smooth over $L_m$ for almost all $m$, by another application
of Proposition~\ref{nonsingular}. This implies that $D_m$ has simple normal crossings for almost all $m$. 
\end{proof}

\section{Limits of $F$-pure thresholds}\label{limits}

The following is our main result. As we will see, it easily implies
the theorem stated in \S 1.

\begin{theorem}\label{main_multiplier}
Let $(L_m)_{m\in\NN}$ be a sequence of fields of positive characteristic such that
$\lim_{m\to\infty}{\rm char}(L_m)=\infty$. We fix a non-principal ultrafilter on $\NN$,
and let $k=[L_m]$. 
If $\fra_m\subseteq L_m[x_1,\ldots,x_n]$ 
are nonzero ideals generated in bounded degree, and if $\fra\subseteq k[x_1,\ldots,x_n]$
is the ideal of polynomials corresponding to $(\fra_m)_{m \geq 1}$, then  for every $\lambda\in\RR_{\geq 0}$ we have
$$\cJ(\fra^{\lambda})_{\rm int}=[\tau(\fra_m^{\lambda})].$$
\end{theorem}

\begin{corollary}\label{cor_main_multiplier}
If $(\fra_m)_{m\in\NN}$ and $\fra$ are as in the above theorem, and $\fra_m$ vanishes at the origin
for almost all $m$, then
$$\lct_0(\fra)={\rm sh}([\fpt_0(\fra_m)]).$$
\end{corollary}

\begin{proof}
Note first that since $\fra_m\subseteq (x_1,\ldots,x_n)L_m[x_1,\ldots,x_n]$ for almost all $m$,
we have $\fra\subseteq (x_1,\ldots,x_n)k[x_1,\ldots,x_n]$. 
By definition, we have 
$$\lct_0(\fra)=\min\{\lambda\in\RR_{\geq 0}\mid \cJ(\fra^{\lambda})\subseteq (x_1,\ldots,x_n)\}.$$
Since $\cJ(\fra^{\lambda})\subseteq (x_1,\ldots,x_n)$ if and only if 
$\cJ(\fra^{\lambda})_{\rm int}\subseteq (x_1,\ldots,x_n)_{\rm int}$, it follows from 
Theorem~\ref{main_multiplier} that this is the case if and only if 
$\tau(\fra_m^{\lambda})\subseteq (x_1,\ldots,x_n)$ for almost all $m$. This is further  equivalent to $\lambda\geq \fpt_0(\fra_m)$ for almost all $m$. We conclude that $\lct_0(\fra)\geq\fpt_0(\fra_m)$ for almost all $m$. In addition, for every $\epsilon\in\RR_{>0}$, 
we have $\cJ(\fra^{\lct_0(\fra) - \epsilon})_{\rm int} \not\subseteq (x_1,\ldots,x_n)_{\rm int}$, and using 
again Theorem~\ref{main_multiplier} we deduce that $\tau( \fra_m^{\lct_0(\fra) - \epsilon}) \not\subseteq(x_1,\ldots,x_n)$ for almost all $m$. By definition, this means that $\fpt_0(\fra_m)\geq\lct_0(\fra)-\epsilon$ for almost all $m$. 
This proves the assertion in the corollary.
\end{proof}

The result stated in \S 1 is an easy consequence of the above corollary.

\begin{proof}[Proof of Theorem~\ref{main_intro}]
Suppose first that we have a sequence $(c_m)_{m\in\NN}$ with $c_m\in {\mathcal Fpt}(p_m)_n$
for all $m$, and such that $\lim_{m\to\infty}p_m=\infty$ and $c=\lim_{m\to\infty} c_m$. We need to show that $c\in {\mathcal Lct}_n$. By Proposition~\ref{dense}, we may assume that 
there are algebraically closed fields $L_m$ of characteristic $p_m$, and ideals
$\fra_m\subseteq L_m[x_1,\ldots,x_n]$ vanishing at the origin, such that $c_m=\fpt_0(\fra_m)$. 
For every $d$, let $\fra_m^{(d)}=\fra_m+(x_1,\ldots,x_n)^d$. 
It follows from Proposition~\ref{prop2_section1} that $|\fpt_0(\fra_m)-\fpt_0(\fra_m^{(d)})|
\leq\frac{n}{d}$.

Let $\cU$ be a non-principal ultrafilter on $\NN$. We put $k=[L_m]$, and for every $d$, we denote
by $\fra^{(d)}\subseteq k[x_1,\ldots,x_n]$ the ideal of polynomials associated to the sequence
of ideals generated in bounded degree $(\fra_m^{(d)})_{m\in\NN}$. 
Given any $\epsilon\in\RR_{>0}$, let $d\gg 0$ be such that $\frac{n}{d}<\epsilon$.
By Corollary~\ref{cor_main_multiplier}, we have
$|\fpt_0(\fra_m^{(d)})-\lct_0(\fra^{(d)})|<\epsilon-\frac{n}{d}$ for almost all $m$.
Therefore $|\fpt_0(\fra_m)-\lct_0(\fra^{(d)})|<\epsilon$ for infinitely many $m$.
Since this holds for every $\epsilon\in\RR_{>0}$, we conclude that $c$ lies in the closure of
$\{\lct_0(\fra^{(d)})\mid d\geq 1\}$. As we have mentioned in \S\ref{subsection_mult_ideals},
${\mathcal Lct}_n$ is closed, hence $c\in {\mathcal Lct}_n$. 

In order to prove the converse, let us consider $c\in {\mathcal Lct}_n$. Consider a sequence of prime integers
$(p_m)_{m\in\NN}$ with limit infinity, and let $L_m$ be an algebraically closed field of characteristic
$p_m$. We fix, as above, a non-principal ultrafilter on $\NN$, and let $k=[L_m]$. As pointed out
in \S\ref{subsection_mult_ideals}, since $k$ is algebraically closed, we can 
find an ideal $\frb\subset k[x_1,\ldots,x_n]$ vanishing at the origin, such that $c=\lct_0(\frb)$. 
Let us write $\frb_{\rm int}=[\frb_m]$.  It follows from Corollary~\ref{cor_main_multiplier}
that $c$ is the limit of a suitable subsequence of $(\fpt_0(\frb_m))_{m\in\NN}$. This completes
the proof of the theorem. Note that the second implication also follows from the results of
\cite{HY} discussed in the Introduction. 
\end{proof}

Before giving the proof of Theorem~\ref{main_multiplier}, we describe the approach from 
\cite{HY} for proving the equality of multiplier ideals with test ideals in a fixed positive characteristic. 
The main ingredients are due independently to Hara \cite{Ha}
and Mehta and Srinivas \cite{MehtaSrinivas}. 
We simplify somewhat the approach in \cite{HY}, 
avoiding the use of local cohomology, which is important in our non-local setting.

Suppose that $L$ is a perfect field of positive characteristic $p$, and $W$ is 
a smooth, irreducible, $n$-dimensional affine scheme over $L$. We consider a nonzero ideal
$\frb$ on $W$, and suppose that we have given a log resolution $\pi\colon \widetilde{W}\to W$
of $\frb$. Let $Z$ be the effective divisor on $\widetilde{W}$ such that 
$\widetilde{\frb}:=\frb\cdot\cO_{\widetilde{W}}=\cO_{\widetilde{W}}(-Z)$, and let
$E=E_1+\ldots+E_N$ be a simple normal crossings divisor on $\widetilde{W}$ such that
both $K_{\widetilde{W}/W}$ and $Z$ are supported on $E$. For every $\lambda\geq 0$, we put
$\cJ(\frb^{\lambda})=\pi_*\cO_{\widetilde{W}}(K_{\widetilde{W}/W}-\lfloor\lambda Z\rfloor)$
(it is irrelevant for us whether this is independent of the given resolution). 
In this setting, it is shown in \cite{HY} that the test ideals are always contained in the multiplier ideals.

\begin{proposition}\label{one_inclusion}
With the above notation, we have
$\tau(\frb^{\lambda})\subseteq\cJ(\frb^{\lambda})$ for all $\lambda\in\RR_{\geq 0}$.
\end{proposition}

\begin{proof}
We give a proof using the description of test ideals at the end of \S 2, since the approach
will be relevant also when considering  the reverse inclusion.
We show that
\begin{equation}\label{eq1_one_inclusion}
(\frb^m)^{[1/p^e]}\subseteq \cJ(\frb^{m/p^e})
\end{equation}
for every $m\geq 0$ and $e\geq 1$. This is enough: given $\lambda\in\RR_{\geq 0}$,
we have for $e\gg 0$
$$\tau(\frb^{\lambda})=(\frb^{\lceil \lambda p^e\rceil})^{[1/p^e]}\subseteq
\cJ(\frb^{\lceil \lambda p^e\rceil/p^e})=\cJ(\frb^{\lambda}).$$
Note that the last equality follows from the fact that $0\leq\frac{\lceil \lambda p^e\rceil}{p^e}-\lambda\ll 1$
for $e\gg 0$.

The commutative diagram (\ref{diag3_1}) induces a commutative diagram 
\begin{equation}\label{diag4_1}
\begin{CD}
F^e_*(\omega_W)@>{t_{W}^e}>>\omega_W\\
@V{\eta=F^e_*(\rho)}VV @VV{\rho}V\\
F^e_*\pi_*(\omega_{\widetilde{W}})@>{\pi_*(t_{\widetilde{W}}^e)}>>\pi_*(\omega_{\widetilde{W}}),
\end{CD}
\end{equation}
where the vertical maps are isomorphisms. Note that $t_{\widetilde{W}}^e$ induces
a (surjective) map $F^e_*(\omega_{\widetilde{W}}(-mZ))\to \omega_{\widetilde{W}}(-
\lfloor \frac{m}{p^e}Z\rfloor)$, and thus a map $F^e_*\pi_*\left(\omega_{\widetilde{W}}(-mZ)\right)\to \pi_*(\omega_{\widetilde{W}}(-
\lfloor \frac{m}{p^e}Z\rfloor))$. Since
$$(\frb^m)^{[1/p^e]}\omega_W=t_W^e(F^e_*(\frb^m\omega_W))$$
and $\eta(F^e_*(\frb^m\omega_W))\subseteq F^e_*\pi_*(\omega_{\widetilde{W}}(-mZ))$,
while 
$$\rho^{-1}(\pi_*(\omega_{\widetilde{W}}(-
\lfloor \frac{m}{p^e}Z\rfloor)))=\cJ(\frb^{m/p^e})\omega_W,$$
we see that (\ref{eq1_one_inclusion}) follows from the fact that
$t^e_W(\frb^m\omega_W)=(\frb^m)^{[1/p^e]}\omega_W$, and  the commutativity of (\ref{diag4_1}).
\end{proof}
 
We now explain a criterion for the reverse inclusion $\cJ(\frb^{\lambda})
\subseteq\tau(\fra^{\lambda})$ to hold. We start with the following proposition.

\begin{proposition}\label{crit1}
Suppose that $\widetilde{W}$ is a smooth, irreducible, $n$-dimensional variety over the perfect field $L$ of positive characteristic $p$. If $E$ is a simple normal crossings divisor on 
$\widetilde{W}$,
and $G$ is a divisor supported on $E$ such that $-G$ is effective, then the canonical morphism
\begin{equation}\label{eq1_crit1}
\Gamma\left(\widetilde{W}, F^e_*\left(\omega_{\widetilde{W}}(\lceil p^eG\rceil)\right)\right)\to
\Gamma\left(\widetilde{W}, \omega_{\widetilde{W}}(\lceil G\rceil)\right)
\end{equation}
is surjective for every $e\geq 1$, provided that
the following two conditions hold:
\begin{enumerate}
\item[(A)] $H^i(\widetilde{W},\Omega_{\widetilde{W}}^{n-i}({\rm log}\,E)(-E+\lceil p^{\ell}G\rceil))=0$ for all $i\geq 1$ and $\ell\geq 1$.
\item[(B)] $H^{i+1}(\widetilde{W},\Omega_{\widetilde{W}}^{n-i}({\rm log}\,E)
(-E+\lceil p^{\ell}G\rceil))=0$ for all $i\geq 1$ and $\ell\geq 0$.
\end{enumerate}
\end{proposition}

This is applied as follows. Suppose that $\lambda\in\RR_{\geq 0}$ is fixed, and we have $\mu>\lambda$
such that $\cJ(\frb^{\lambda})=\cJ(\frb^{\mu})$ (note that if $Z=\sum_ia_iE_i$, then
it is enough to take $\mu$ such that $\mu<\frac{\lfloor\lambda a_i\rfloor+1-\lambda a_i}{a_i}$
for all $i$ with $a_i>0$). Let us consider now 
a $\QQ$-divisor
$D$ on $\widetilde{W}$ such that $D$ is ample over $W$, and $-D$
is effective\footnote{Such a divisor always exists: if we express $\widetilde{W}$
as the blow-up of $W$ along a suitable ideal, then we may take $D$ to be the the negative of the
exceptional divisor.}. We will apply the above proposition with $G=\mu(D-Z)$.
We may and will assume that
$\lceil G\rceil=\lceil -\mu Z\rceil$ (again this condition only depends on $\mu$
and the coefficients of $Z$; since $-D$ is effective, it is always satisfied if we replace $D$ by $\epsilon D$, with
$0<\epsilon\ll 1$).

\begin{proposition}\label{crit2}
With the above notation, if (\ref{eq1_crit1}) is surjective for every $e\geq 1$, then 
$\cJ(\frb^{\lambda})\subseteq\tau(\frb^{\lambda})$.
\end{proposition}

\begin{proof}
We use again the commutative diagram (\ref{diag4_1}). This induces 
a commutative diagram
\begin{equation}\label{diag4_2}
\begin{CD}
F^e_*\pi_*(\omega_{\widetilde{W}}(\lceil p^eG\rceil)) @>{\pi_*(t_{\widetilde{W}}^e)}>>
\pi_*(\omega_W(\lceil G\rceil))=\pi_*(\omega_W(-\lfloor\mu Z\rfloor))\\
@V{\eta^{-1}}VV @VV{\rho^{-1}}V\\
F^e_* \omega_W@>{t_W^e}>>\omega_W
\end{CD}
\end{equation}
in which the top horizontal map is surjective by assumption (recall that $W$ is affine), and the image of 
the right vertical map is $\cJ(\frb^{\mu})\omega_W$. The image of the left vertical map
can be written as $F_*^e(J_e\omega_W)$, where $J_e=
\pi_*\cO_{\widetilde{W}}(K_{\widetilde{W}/W}+\lceil p^eG\rceil)
$, and we deduce from the commutativity of (\ref{diag4_2}) that $\cJ(\frb^{\mu})\subseteq J_e^{[1/p^e]}$. 
By Lemma~\ref{lem_multiplier} below, there is $r$ such that
$\cJ(\frb^m)\subseteq\frb^{m-r}$ for every $m\geq r$. 
Since $-D$ is effective, by letting $e\gg 0$, we get
$$J_e=\pi_*\cO_{\widetilde{W}}(K_{\widetilde{W}/W}-\lfloor p^e\mu(Z-D)\rfloor)
\subseteq \cJ(\frb^{\mu p^e})\subseteq\frb^{\lfloor \mu p^e\rfloor -r},$$
and therefore
$$\cJ(\frb^{\lambda})=\cJ(\frb^{\mu})\subseteq 
(\frb^{\lfloor \mu p^e\rfloor -r})^{[1/p^e]}\subseteq\tau(\frb^{\frac{\lfloor\mu p^e\rfloor-r}{p^e}})
\subseteq\tau(\fra^{\lambda}),$$
since $\lim_{e\to\infty}\frac{\lfloor\mu p^e\rfloor-r}{p^e}=\mu>\lambda$. This completes the proof of the proposition.
\end{proof}

\begin{lemma}\label{lem_multiplier}
With the above notation, there is $r$ such that $\cJ(\frb^m)\subseteq\frb^{m-r}$
for every integer $m\geq r$. 
\end{lemma}

\begin{proof}
It is enough to prove, more generally, that for every coherent sheaf $\cF$ on $\widetilde{W}$,
the graded module $M:=\oplus_{m\geq 0}\Gamma(\widetilde{W},\cF(-mZ))$
is finitely generated over the Rees algebra $S:=\oplus_{m\geq 0}\frb^m$. 
We may factor $\pi$ as $\widetilde{W}\overset{g}\to B\overset{f}\to W$, where $B$
is the normalized blow-up of $W$ along $\frb$ (that is, $B={\rm Proj}(S')$, where 
$S'$ is the normalization of $S$). The line bundle $\frb\cdot\cO_B=\cO_B(-T)$
is ample over $W$, and using the projection formula we see that
$M=\oplus_{m\geq 0}\Gamma(B,\pi_*(\cF)\otimes\cO_B(-mT))$
is finitely generated over $S'=\oplus_{m\geq 0}\Gamma(B,\cO_B(-mT))$.
Since $S'$ is a finite $S$-algebra, it follows that $M$ is a finitely generated $S$-module.
\end{proof}

We recall, for completeness, the proof of Proposition~\ref{crit1}, which makes use of the de Rham complex $\Omega_{\widetilde{W}}^{\bullet}({\rm log}(E))$ with log poles along the simple normal crossings divisor $E$. Note that while this complex does not have
$\cO_{\widetilde{W}}$-linear differentials, its Frobenius push-forward
$F_*\Omega_{\widetilde{W}}^{\bullet}({\rm log}(E))$ does have this property. 
In particular, we may tensor this complex with line bundles. If $\cL$
is a line bundle, then by the projection formula we have
$$(F_*\Omega_{\widetilde{W}}^i({\rm log}\,E))\otimes\cL\simeq F_*(\Omega_{\widetilde{W}}^i
({\rm log}\,E)\otimes\cL^p).$$
The following facts are the key ingredients in the proof of Proposition~\ref{crit1}.

\begin{enumerate}
\item The Cartier isomorphism: there is a canonical isomorphism (see \cite[Theorem 1.2]{DI})
\[ C^{-1}\colon \Omega^i_{\widetilde{W}}(\log E) \simeq \cH^i 
F_* (\Omega^\bullet_{\widetilde{W}}(\log E)). \]
\item Insensitivity to small effective  twists: 
suppose that $B$ is an effective divisor supported on $E$, with all coefficients less than $p$.
We have a twisted de Rham complex with log poles $\Omega_{\widetilde{W}}^{\bullet}({\rm log}\,E)(B)$ (it is enough to check that the differential of the de Rham complex of meromorphic 
differential forms on $X$ preserves these subsheaves). In this case, the natural inclusion
 $\Omega^\bullet_{\widetilde{W}}({\rm log}\, E) \hookrightarrow \Omega_{\widetilde{W}}^{\bullet}({\rm log}\,E)(B)$ is a quasi-isomorphism (see \cite[Lemma~3.3]{Ha} or 
 \cite[Corollary~4.2]{MehtaSrinivas} for a proof). Combining this with the Cartier isomorphism, we find
 \begin{equation}\label{eq3_crit1}
 \Omega^i_{\widetilde{W}}({\rm log}\, E) \simeq \cH^i \left(F_*\left(\Omega_{\widetilde{W}}^{\bullet}({\rm log}\,E)(B)\right)\right). 
 \end{equation}
\end{enumerate}

\begin{proof}[Proof of Proposition~\ref{crit1}]
Note first that it is enough to prove the case $e=1$. Indeed, if $\alpha_{G,e}$ is the morphism
(\ref{eq1_crit1}), we see that $\alpha_{G,e}=\alpha_{G,1}\circ\alpha_{pG,1}\circ\ldots
\alpha_{p^{e-1}G,1}$, and the hypothesis implies that we may apply the condition for $e=1$
to each of $G,pG,\ldots,p^{e-1}G$. Therefore from now on we assume that $e=1$
(and in this case we will only need condition (A) for $\ell=1$ and condition (B) for
$\ell=0$).

Let $B:=(p-1)E+\lceil pG\rceil-p\lceil G\rceil=(p-1)E+p\lfloor -G\rfloor-\lfloor -pG\rfloor$.
Since $-G$ is effective, it follows from the second expression that $B$ is effective, and its coefficients are less than $p$. Let $K^{\bullet}:=F_*\Omega_{\widetilde{W}}^{\bullet}({\rm log}\,E)(-E+\lceil pG\rceil)$. By tensoring (\ref{eq3_crit1}) with $\cO_{\widetilde{W}}(-E+\lceil G\rceil)$,
and using the projection formula, we get
$$\Omega_{\widetilde{W}}^i({\rm log}\,E)(-E+\lceil G\rceil)
\simeq\cH^i\left(F_*\left(\Omega_{\widetilde{W}}^{\bullet}({\rm log}\,E)(B-pE+p\lceil G\rceil)
\right)\right)=\cH^i(K^{\bullet}).$$
Note that the morphism $\alpha_{G,1}$ is identified to $\Gamma(\widetilde{W},K^n)
\to \Gamma(\widetilde{W}, \cH^n(K^{\bullet}))$. It is then straightforward to show, by breaking 
$K^{\bullet}$ into short exact sequences, and using the corresponding long exact sequences for cohomology, that $\alpha_{G,1}$ is surjective if
$H^i(\widetilde{W},K^{n-i})=0$ and $H^{i+1}(\widetilde{W}, \cH^{n-i}(K^{\bullet}))=0$
for all $i\geq 1$. By what we have seen, these are precisely conditions (A) with $\ell=1$ and (B)
with $\ell=0$.
\end{proof}

We will also make use of the following version of the Kodaira-Akizuki-Nakano
vanishing theorem (in characteristic zero).

\begin{theorem}\label{vanishing_theorem}
Let $Y$ be a smooth, irreducible variety over a field $k$ of characteristic zero. 
If $Y$ is projective over an affine scheme $X$, 
$E$ is a reduced simple normal crossings divisor on $Y$, 
and $G$ is a $\QQ$-divisor on $Y$ such that $G-\lfloor G\rfloor$ is supported on $E$ and $G$ is ample over $X$, then 
$$H^i\left(Y, \Omega^j_Y({\rm log}\,E)(-E+\lceil G\rceil)\right)=0\,\,\text{if}\,\,i+j>\dim(X).$$
\end{theorem}

\begin{proof}
This is proved when ${\rm char}(k)=p>0$ in \cite[Corollary~3.8]{Ha} under the assumption
that $p>\dim(X)$ and that both $Y$ and $E$ admit a lifting to the second ring of Witt vectors
$W_2(k)$ of $k$. The proof relies on an application of the results from \cite{DI}.
It is then standard to deduce the assertion in characteristic zero (see, for example, the proof of
\cite[Corollary~2.7]{DI}).
\end{proof}

We can now give the proof of our main result.

\begin{proof}[Proof of Theorem~\ref{main_multiplier}]
Let $p_m={\rm char}(L_m)$. We have by hypothesis $\lim_{m\to\infty}p_m=\infty$, hence
${\rm char}(k)=0$. In particular, there is a log resolution $\pi\colon Y\to X=\AAA_k^n$
of $\fra$. We write $\fra\cdot\cO_Y=\cO_Y(-Z)$, and let $E$ be a simple normal crossings divisor on $Y$ such that both $Z$ and $K_{Y/X}$ are supported on $E$.
Let $[\pi_m]\colon [Y_m]\to [X_m]=[\AAA_{L_m}^n]$ be the corresponding morphism
of internal schemes. It follows from Proposition~\ref{various_properties2} (iii) that
if $Z_{\rm int}=[Z_m]$, then $\fra_m\cdot\cO_{Y_m}=\cO(-Z_m)$ for almost all $m$.
On the other hand,
it is easy to deduce from Proposition~\ref{nonsingular} that
$(K_{Y/X})_{\rm int}=[K_{Y_m/X_m}]$. If $E_{\rm int}=[E_m]$, then $E_m$ has simple normal crossings for almost all $m$ by Proposition~\ref{SNC}, and we conclude that $\pi_m$ is a log resolution of $\fra_m$ for almost all $m$. Moreover, if we use $\pi_m$ to define $\cJ(\fra_m^{\lambda})$ on $X_m$, then
we have $\cJ(\fra^{\lambda})_{\rm int}=[\cJ(\fra_m^{\lambda})]$ by 
Proposition~\ref{various_properties2} (iv). 

For every $m$ such that $\pi_m$ gives a log resolution of $\fra_m$ we have
$\tau(\fra_m^{\lambda})\subseteq\cJ(\fra_m^{\lambda})$ by Proposition~\ref{one_inclusion}.
We now choose a rational number $\mu>\lambda$ such that $\cJ(\fra^{\lambda})=\cJ(\fra^{\mu})$,
so that $\cJ(\fra_m^{\lambda})=\cJ(\fra_m^{\mu})$ for almost all $m$. 
We also choose a $\QQ$-divisor $D$ supported on $E$ such that $-D$ is effective, $D$ is ample over $X$, and $\lceil \mu(D-Z)\rceil=\lceil -\mu Z\rceil$. We write $G=\mu(D-Z)$,
and denote by $D_m$ and respectively $G_m$ the corresponding divisors on $X_m$.
It is clear that for almost all $m$ the divisor $-D_m$ is effective, $D_m$ is ample over $X_m$
(see Proposition~\ref{general_properties} (v)), and  $\lceil G_m\rceil=\lceil -\mu Z_m\rceil$.
We deduce from
Propositions~\ref{crit1} and 
\ref{crit2} that $\cJ(\fra_m^{\lambda})\subseteq\tau(\fra_m^{\lambda})$
if the following conditions hold:
\begin{enumerate}
\item[(${\rm A}_m$)] $H^i(Y_m,\Omega_{Y_m}^{n-i}({\rm log}\,E_m)(-E_m+\lceil p_m^{\ell}G_m\rceil))=0$ for all $i\geq 1$ and $\ell\geq 1$.
\item[(${\rm B}_m$)] $H^{i+1}(Y_m,\Omega_{Y_m}^{n-i}({\rm log}\,E_m)
(-E_m+\lceil p_m^{\ell}G_m\rceil))=0$ for all $i\geq 1$ and $\ell\geq 0$.
\end{enumerate}
It follows that in order to complete the proof, it is enough to show that conditions 
(${\rm A}_m$) and (${\rm B}_m$) hold for almost all $m$.

Note first that by Theorem~\ref{vanishing_theorem}, we have 
$H^{i+1}\left(Y,\Omega^{n-i}_Y({\rm log}\,E)(-E+\lceil G\rceil)\right)=0$ for all $i\geq 0$.
Using Proposition~\ref{various_properties2} (iv), we deduce that
$H^{i+1}\left(Y_m,\Omega^{n-i}_{Y_m}({\rm log}\,E_m)(-E_m+\lceil G_m\rceil)\right)=0$ 
for all $i\geq 0$ and almost all $m$ (since these groups vanish automatically
when $i\geq n$, we only need to consider finitely many such $i$.
This takes care of the condition (${\rm B}_m$) for $\ell=0$. 

We now treat the remaining conditions. Let us fix  a positive integer $d$ such that 
$dG$ is an integral divisor.  Let $\cF_1,\ldots,\cF_M$ denote the sheaves
$\Omega^i_Y({\rm log}\,E)(-E+\lceil sG\rceil)$, for integers $0\leq s\leq n$
and $0\leq j\leq d-1$. Since $dG$ is ample over the affine variety $X$, there is 
$j_0$ such that $H^i(Y,\cF_t(jdG))=0$ for every $j\geq j_0$, every $i\geq 1$ and every 
$t\leq M$.
If $m$ is such that $p_m\geq (j_0+1)d$, and if for $\ell\geq 1$ we take  $s$ with
$0\leq s\leq d-1$ such that $p_m^{\ell}\equiv s$ (mod $d$), then
$$\lceil p_m^{\ell}G_m\rceil=\frac{p_m^{\ell}-s}{d}(dG)+\lceil sG\rceil\,\,\text{and}\,\,
\frac{p_m^{\ell}-s}{d}\geq \frac{p_m-s}{d}\geq j_0.$$
We deduce from Corollary~\ref{Serre_vanishing} that the vanishings in 
(${\rm A}_m$) and (${\rm B}_m$) hold when $\ell\geq 1$ for almost all $m$
(note that for such $m$ we may assume that $p_m\geq (j_0+1)d$). This completes the
proof of the theorem. 
\end{proof}

\providecommand{\bysame}{\leavevmode \hbox \o3em
{\hrulefill}\thinspace}

\end{document}